\documentclass{article}[12pt]

\usepackage{graphicx}
\usepackage{mathptmx}
\usepackage{amssymb}
\usepackage{latexsym} 
\usepackage{amsmath}
\usepackage{amssymb}
\usepackage{graphics}
\usepackage{marvosym}
\usepackage{caption}

\begin{document}

\vspace{0.5cm}

\centerline{\bf\Large A Short Note on  Contracting Self-Similar Solutions}
\centerline{\bf\Large of the Curve Shortening Flow}

\vspace{0.5cm}

\centerline{Lucas Z. Veeravalli, Emma H. Veeravalli, Alain R. Veeravalli}

\vspace{0.5cm}

\vspace{0.5cm}

\centerline{Universit\'e d'Evry-Val d'Essonne, D\'epartement de Math\'ematiques}
\centerline{23 Boulevard de France, 91037 Evry Cedex, France}
\centerline{\Letter\ Alain.Veeravalli@univ-evry.fr}

\vspace{0.5cm}

\centerline{Received: 30 June 2015}

\vspace{0.5cm}

\centerline{\bf\large Abstract}
\noindent By the curve shortening flow, the only closed contracting self-similar solutions are circles: we give a very short and intuitive geometric proof of this basic and classical result using an idea of Gage \cite{ga}.\newline


\noindent MSC 53C44, 53A04.

\section{Introduction}
\label{intro}

\noindent Let ${\boldsymbol\gamma}$ be a smooth closed curve parametrized by arc length, embedded in the Euclidean plane endowed with its canonical inner product denoted by a single point. A one-parameter smooth family of plane closed curves $\left({\bf{\boldsymbol\gamma}}(\cdot,t)\right)_t$ with initial condition ${\bf{\boldsymbol\gamma}}(\cdot,0)={\bf{\boldsymbol\gamma}}$ is said to evolve by the curve shortening flow (CSF for short) if

\begin{equation}\label{csf1}
\frac{\partial{\bf{\boldsymbol\gamma}}}{\partial t}=\kappa {\bf n}
\end{equation}

\noindent where $\kappa$ is the signed curvature and ${\bf n}$ the inward pointing unit normal. By the works of Gage, Hamilton and Grayson any embedded closed curve evolves to a convex curve (or remains convex if so) and shrinks to a point in finite time\footnote[1]{The reader could find a dynamic illustration of this result on the internet page http://a.carapetis.com/csf/}.

\noindent In this note, we are interested by {\it self-similar} solutions that is solutions which shapes change homothetically during the evolution. This condition is equivalent to say, after a suitable parametrization, that

\begin{equation}\label{csf2}
\kappa=\varepsilon{\bf{\boldsymbol\gamma}}\cdot{\bf n}
\end{equation}

\noindent with $\varepsilon=\pm 1$. If $\varepsilon=-1$ (resp. $+1$), the self-similar family is called {\it contracting} (resp. {\it expanding}). For instance, for any positive constant $C$ the concentric circles $\left(s\mapsto \sqrt{1-2t}\,(\cos s,\sin s)\right)_t$ form a self-similar contracting solution of the CSF shrinking to a point in finite time and as a matter of fact, there is no more example than this one:\newline

{\it by the curve shortening flow, the only closed embedded contracting self-similar solutions are circles\footnote[2]{The nonembedded closed curves were studied and classified by Abresch and Langer\cite{al}}.}\newline

\noindent To the author knowledge, the shortest proof of this was given by Chou-Zhu \cite{cz} by evaluating a clever integral. The proof given here is purely geometric and based on an genuine trick used by Gage in \cite{ga}.

\section{A geometric proof}
\label{sec:1}

\noindent Let ${\bf{\boldsymbol\gamma}}$ be a closed, simple embedded plane curve, parametrized by arclength $s$, with signed curvature $\kappa$. By reversing the orientation if necessary, we can assume that the curve is counter-clockwise oriented. The length of ${\bf{\boldsymbol\gamma}}$ is denoted by $L$, the compact domain enclosed by ${\bf{\boldsymbol\gamma}}$ will be denoted by $\Omega$ with area $A$ and the associated moving Frenet frame by $({\bf t},{\bf n})$. Let ${\boldsymbol\gamma_t}={\boldsymbol\gamma}(t,\cdot)$ be the one parameter smooth family solution of the CSF, with the initial condition ${\bf{\boldsymbol\gamma}}_0={\bf{\boldsymbol\gamma}}$.\newline

\noindent Multiplying $(1)$ by ${\bf n}$, we obtain

\begin{equation}\label{csf2}
\frac{\partial{\bf{\boldsymbol\gamma}}}{\partial t}\cdot{\bf n}=\kappa
\end{equation}

\noindent Equations $(\ref{csf1})$ and $(\ref{csf2})$ are equivalent: from $(\ref{csf2})$, one can look at a reparametrization $(t,s)\mapsto \varphi(t,s)$ such that $\widetilde{\bf{\boldsymbol\gamma}}(t,s)={\bf{\boldsymbol\gamma}}(t,\varphi(t,s))$ satisfies $(1)$. A simple calculation leads to an ode on $\varphi$ which existence is therefore guaranteed \cite{cz}. From now, we will deal with equation $(\ref{csf2})$.\newline

\noindent If a solution $\bf{\boldsymbol\gamma}$ of $(\ref{csf2})$ is self-similar, then there exists a non-vanishing smooth function $t\mapsto \lambda(t)$ such that ${\boldsymbol\gamma}_t(s)=\lambda(t)\,{\boldsymbol\gamma}(s)$. By $(2)$, this leads to $\lambda'(t)\,{\boldsymbol\gamma}(s)\cdot{\bf n}({\boldsymbol\gamma}_t(s))=\kappa({\boldsymbol\gamma}_t(s))$, that is $\lambda'(t)\lambda(t)\,{\boldsymbol\gamma}(s)\cdot{\bf n}(s)=\kappa(s)$. The function $s\mapsto{\boldsymbol\gamma}(s)\cdot{\bf n}(s)$ must be non zero at some point (otherwise $\kappa$ would vanish everywhere and ${\boldsymbol\gamma}$ would be a line) so the function $\lambda'\lambda$ is constant equal to a real $\varepsilon$ which can not be zero. By considering the new curve $s\mapsto \sqrt{\vert \varepsilon\vert}{\boldsymbol\gamma}\left(s/\sqrt{\vert \varepsilon\vert}\right)$ which is still parametrized by arc length, we can assume that $\varepsilon=\pm 1$. In the sequel we will assume that ${\boldsymbol\gamma}$ is contracting, that is $\varepsilon=-1$ which says that we have the fundamental relation:

\begin{equation}\label{fund-rela}
\kappa+{\boldsymbol\gamma}\cdot {\bf n}=0
\end{equation}

\noindent An immediate consequence is the value of $A$: indeed, by the divergence theorem and the turning tangent theorem,

$$A=\frac{1}{2}\int_{{\boldsymbol\gamma}} \left(x\,dy-y\,dx\right)=-\frac{1}{2}\int_0^L {\boldsymbol\gamma}(s)\cdot{\bf n}(s)\,ds=\frac{1}{2}\int_0^L k(s)\,ds=\pi$$

\noindent Therefore, our aim will be to prove that $L=2\pi$ and we will conclude by using the equality case in the isoperimetric inequality.\newline

\noindent The second remark is that the curve is an oval or strictly convex: indeed, by differentiating $(\ref{fund-rela})$ and using Frenet formulae, we obtain that $\kappa'=\kappa\,{\boldsymbol\gamma}\cdot{\boldsymbol\gamma}'$ which implies that $\kappa=Ce^{\vert{\bf{\boldsymbol\gamma}}\vert^2/2}$ for some non zero constant $C$. As the rotation index is $+1$, $C$ is positive and so is $\kappa$.\newline

\subsection{Polar tangential coordinates}
\label{sec:2}

\noindent As equation $(\ref{fund-rela})$ is invariant under Euclidean motions, we can assume that the origin $O$ of the Euclidean frame lies within $\Omega$ with axis $[Ox)$ meeting ${\boldsymbol\gamma}$ orthogonally. We introduce the angle function $\theta$ formed by $-{\bf n}$ with the $x$-axis as shown in the figure below:


\begin{figure}[h]
\includegraphics*[scale=0.8,bb=120 350 550 750]{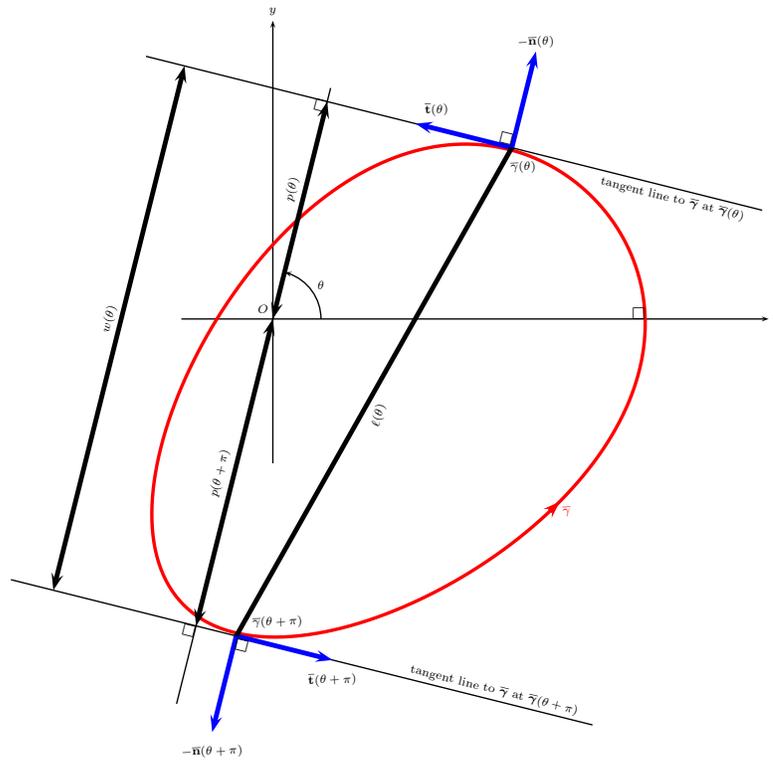}
\caption{Polar tangential coordinates}
\label{fig:1}
\end{figure}

\noindent Then $-{\bf n}(s)=(\cos\theta(s),\sin\theta(s))$. Differentiating this equality, we obtain $\kappa(s)\,{\bf t}(s)=\theta'(s)\,{\bf t}(s)$, that is $\theta(s)=\int_0^s \kappa(u)\,du$ since $\theta(0)=0$. As $\theta'=\kappa\geqslant C>0$, $\theta$ is a strictly increasing function on $\mathbb{R}$ onto $\mathbb{R}$. So $\theta$ can be choosen as a new parameter and we set $\overline{{\bf{\boldsymbol\gamma}}}(\theta)=(\overline{x}(\theta),\overline{y}(\theta))={\bf{\boldsymbol\gamma}}(s)$, $\overline{\bf t}(\theta)=(-\sin\theta,\cos\theta)$, $\overline{\bf n}(\theta)=(-\cos\theta,-\sin\theta)$ and we consider the function $p$ defined by $p(\theta)=-\overline{{\bf{\boldsymbol\gamma}}}(\theta)\cdot \overline{{\bf n}}(\theta)$. As $\theta(s+L)=\theta(s)+2\pi$, we note that $\overline{\bf t}$, $\overline{\bf n}$ and $p$ are $2\pi$-periodic functions. The curve $\overline{{\bf{\boldsymbol\gamma}}}$ is regular but not necessarily parametrized by arc length because $\overline{{\bf{\boldsymbol\gamma}}}'(\theta)=\frac{1}{k(s)}{\bf{\boldsymbol\gamma}}'(s)$ and we note $\overline{\kappa}$ its curvature. By definition, we have

\begin{equation}\label{p=}
\overline{x}(\theta)\cos\theta+\overline{y}(\theta)\sin\theta=p(\theta)
\end{equation}

\noindent which, by differentiation w.r.t. $\theta$, gives

\begin{equation}\label{p'=}
-\overline{x}(\theta)\sin\theta+\overline{y}(\theta)\cos\theta=p'(\theta)
\end{equation}

\noindent Thus,

\begin{equation}\label{xy}
\left\{
    \begin{array}{ll}
      \overline{x}(\theta) &= p(\theta)\cos\theta-p'(\theta)\sin\theta\\
      \overline{y}(\theta) &= p(\theta)\sin\theta+p'(\theta)\cos\theta\\
    \end{array}
  \right.
\end{equation}

\noindent Differentiating once more, we obtain

\begin{equation}\label{x'y'}
\left\{
    \begin{array}{ll}
      \overline{x}'(\theta) &= -\left[p(\theta)+p''(\theta)\right]\sin\theta\\
      \overline{y}'(\theta) &= \hspace{0.32cm}\left[p(\theta)+p''(\theta)\right]\cos\theta\\
    \end{array}
  \right.
\end{equation}

\noindent Since $\gamma$ is counter-clockwise oriented, we have $p+p''>0$.\newline

\noindent Coordinates $(\theta,p(\theta))_{0\leqslant \theta\leqslant 2\pi}$ are called {\it polar tangential coordinates} and $p$ is the {\it Minkowski support function}. By (\ref{x'y'}), we remark that the tangent vectors at $\overline{{\bf{\boldsymbol\gamma}}}(\theta)$ and $\overline{{\bf{\boldsymbol\gamma}}}(\theta+\pi)$ are parallel. We will introduce the {\it width function} $w$ defined by

$$w(\theta)=p(\theta)+p(\theta+\pi)$$

\noindent which is the distance between the parallel tangent lines at $\overline{{\bf{\boldsymbol\gamma}}}(\theta)$ and $\overline{{\bf{\boldsymbol\gamma}}}(\theta+\pi)$ and we denote by $\ell(\theta)$ the segment joining $\overline{{\bf{\boldsymbol\gamma}}}(\theta)$ and $\overline{{\bf{\boldsymbol\gamma}}}(\theta+\pi)$.\newline

\noindent With these coordinates, the perimeter has a nice expression:

\begin{equation}
L=\int_0^{2\pi}({\overline{x}'^2+\overline{y}'^2})^{1/2}d\theta=\int_0^{2\pi} \left(p+p''\right)d\theta=\int_0^{2\pi} p\,d\theta\tag{\mbox{\rm Cauchy formula}}
\end{equation}

\noindent The curvature $\overline{\kappa}$ of $\overline{{\bf{\boldsymbol\gamma}}}$ is

\begin{equation*}
{\overline\kappa}=\frac{\overline{x}'\overline{y}''-\overline{x}''\overline{y}'}{\left(\overline{x}'^2+\overline{y}'^2\right)^{3/2}}=\frac{1}{p+p''}
\end{equation*}

\noindent and equation $(\ref{fund-rela})$ reads ${\overline\kappa}=p$. So, finally,

\begin{equation}\label{k=p}
{\overline\kappa}=p=\frac{1}{p+p''}
\end{equation}

\subsection{Bonnesen inequality}
\label{sec:3}

\begin{figure}[h]
\centerline{\includegraphics[scale=0.7,bb=-300 0 550 140]{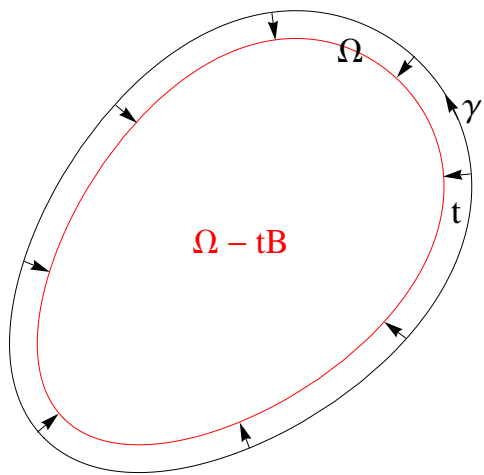}}
\caption{The domain $\Omega-tB$ with positive $t$}
\label{fig:2}
\end{figure}

\noindent If $B$ is the unit ball of the Euclidean plane, it is a classical fact that the area of $\Omega-tB$ (figure \ref{fig:2}) is $A_\Omega(t)=A-Lt+\pi t^2$ \cite{bz}.

\vspace{3.0cm}

\noindent The roots $t_1,t_2$ (with $t_1\leqslant t_2$) of $A_\Omega(t)$ are real by the isoperimetric inequality and they have a geometric meaning: indeed, if $R$ is the circumradius of $\Omega$, that is the radius of the circumscribed circle, and if $r$ is the inradius of $\Omega$, that is the radius of the inscribed circle, Bonnesen \cite{bz,zc} proved in the 1920's a series of inequalities, one of them being the following one:

$$t_1\leqslant r\leqslant R\leqslant t_2$$

\noindent Moreover, and this is a key point in the proof, any equality holds if and only if $\overline{{\bf{\boldsymbol\gamma}}}$ is a circle. We also note that $A_\Omega(t)<0$ for any $t\in(t_1,t_2)$.

\subsection{End of proof}
\label{sec:4}

\noindent {\bf Special case: $\overline{{\bf{\boldsymbol\gamma}}}$ is symmetric w.r.t. the origin $O$}, that is $\overline{{\bf{\boldsymbol\gamma}}}(\theta+\pi)=-\overline{{\bf{\boldsymbol\gamma}}}(\theta)$ for all $\theta\in[0,2\pi]$, which also means that $p(\theta+\pi)=p(\theta)$ for all $\theta\in[0,2\pi]$. So the width function $w$ is twice the support function $p$. As $2r\leqslant w\leqslant 2R$, we deduce that for all $\theta$, $r\leqslant p(\theta)\leqslant R$. If $\overline{{\bf{\boldsymbol\gamma}}}$ is not a circle, then one would derive from Bonnesen inequality that $t_1<r\leqslant p(\theta)\leqslant R<t_2$. So $A_\Omega(p(\theta))<0$ for all $\theta$, that is $\pi p^2(\theta)<Lp(\theta)-\pi$. Multiplying this inequality by $1/p=p+p''$ ($>0$) and integrating on $[0,2\pi]$, we would obtain $\pi L<\pi L$ by Cauchy formula ! By this way, we proved that any symmetric smooth closed curve satisfying $(\ref{fund-rela})$ is a circle. As the area is $\pi$, the length is $2\pi$ of course.

\vspace{0.4cm}

\noindent {\bf General case:} using a genuine trick introduced by Gage \cite{ga}, we assert that\newline

{\it for any oval enclosing a domain of area $A$, there is a segment $\ell(\theta_0)$ dividing $\Omega$ into two subdomains of equal area $A/2$.}\newline

\noindent Proof: let $\sigma(\theta)$ be the area of the subdomain of $\Omega$, bounded by $\overline{{\bf{\boldsymbol\gamma}}}([\theta,\theta+\pi])$ and the segment $\ell(\theta)$. We observe that $\sigma(\theta)+\sigma(\theta+\pi)=A$. We can assume without lost of generality that $\sigma(0)\leqslant A/2$. Then we must have $\sigma(\pi)\geqslant A/2$, and by continuity of $\sigma$ and the intermediate value theorem, there exists $\theta_0$ such that $\sigma(\theta_0)=A/2$ and the segment $\ell(\theta_0)$ proves the claim.\hfill$\Box$

\vspace{0.4cm}

\begin{figure}[h]
\centerline{\includegraphics*[scale=0.8,bb=120 400 550 900]{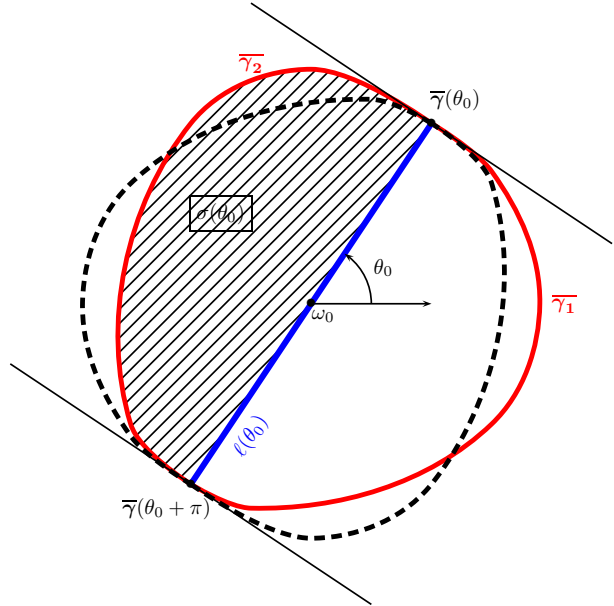}}
\caption{Symmetrization of the curve}
\label{fig:3}
\end{figure}

\noindent Let $\omega_0$ be the center of $\ell(\theta_0)$. If $\overline{{\bf{\boldsymbol\gamma}}}_1$ and $\overline{{\bf{\boldsymbol\gamma}}}_2$ are the two arcs of $\overline{{\bf{\boldsymbol\gamma}}}$ separated by $\ell(\theta_0)$, we denote by $\overline{{\bf{\boldsymbol\gamma}}}_i^s$ ($i=1,2$) the closed curve formed by $\overline{{\bf{\boldsymbol\gamma}}}_i$ and its reflection trough $\omega_0$. Each $\overline{{\bf{\boldsymbol\gamma}}}_i^s$ is a symmetric closed curve and as $\ell(\theta_0)$ joins points of the curve where the tangent vectors are parallel, each one is strictly convex and smooth (figure \ref{fig:3}).

\vspace{-6.0cm}

\noindent Moreover, each $\overline{{\bf{\boldsymbol\gamma}}}_i^s$ satisfies equation $(\ref{fund-rela})$ and encloses a domain of area $2\times A/2=A$. So we can apply the previous case to these both curves and this gives that $\mbox{\rm length}(\overline{{\bf{\boldsymbol\gamma}}}_i^s)=2\pi$ for $i=1,2$, that is $\mbox{\rm length}(\overline{{\bf{\boldsymbol\gamma}}}_i)=\pi$ which in turn implies that $L=2\pi$. So $\overline{{\bf{\boldsymbol\gamma}}}$ (that is ${\bf{\boldsymbol\gamma}}$) is a circle. This proves the theorem.\hfill$\Box$

\end{document}